\documentclass{amsart}
\title{Vanishing Results for Toric Varieties Associated to $GL_n$ and $G_2$}

\usepackage{epic,eepic, graphicx, amssymb}

\newtheorem{definition}{Definition}[section]

\newtheorem{proposition}[definition]{Proposition}
\newtheorem{lemma}[definition]{Lemma}
\newtheorem{remark}[definition]{Remark}

\newcommand{\R}{\mathbb{R}}
\newcommand{\Z}{\mathbb{Z}}
\newcommand{\C}{\mathbb{C}}
\newcommand{\noadd}[1]{}

\input xy

\xyoption{all}

\newcommand{\Dim}{\operatorname{dim}}

\newcommand{\Max}{\operatorname{Max}}

\newcommand{\Conv}{\operatorname{Conv}}

\author{Q\"endrim R. Gashi}
\begin{document}

\begin{abstract}
Toric varieties associated with root systems appeared very naturally in the theory of group compactifications. Here they are considered in a very different context.

We prove the vanishing of higher cohomology groups for certain line bundles on toric varieties associated to $GL_n$ and $G_2$. This can be considered of general interest and it improves the previously known results for these varieties. We also show how these results give a simple proof of a converse to Mazur's inequality for $GL_n$ and $G_2$ respectively. It is known that the latter imply the non-emptiness of some affine Deligne-Lusztig varieties.
 
\end{abstract}
\maketitle
\begin{center}
Dedicated to Scarlett MccGwire and Dr. Christian Duhamel
\end{center}

\section*{Introduction}

Fix an isocrystal and the corresponding Newton vector. Mazur's inequality (see ~\cite{mazur}, pg.658, and the group-theoretic generalization in ~\cite{rapritch}, Theorem 4.2) states that, given a lattice in our isocrystal, the Hodge vector associated to it is greater than or equal to the Newton vector. Conversely, any vector satisfying Mazur's inequality is the Hodge vector for some lattice (see e.g. ~\cite{rapoprt}, Proposition 4.2, and ~\cite{frapo}, Theorem 2). The last result can be regarded as a statement for $GL_n$, but, in ~\cite{krapo}, \S 4, R. Kottwitz and M. Rapoport formulated an analogous version for other groups and they proved it for $GSp_{2n}$ and $GL_n$. They also formulated a combinatorial statement involving respective root systems that would imply the converse to Mazur's inequality in each case (see also sections 4.3 and 4.4 in ~\cite{kot}, for a detailed version). C. Lucarelli in ~\cite{cathy}, Theorem 0.2, showed that this combinatorial statement is true (see Proposition \ref{cathy} below), and therefore a converse to Mazur's inequality was deduced, for all split classical groups. We prove this combinatorial statement for $G_2$ (Theorem A) and hence deduce a converse to Mazur's inequality for $G_2$. We also generalize the combinatorial result for $GL_n$ (Theorem B), which in particular gives a new and simple proof of the converse to Mazur's inequality for $GL_n$.

The results obtained here are proved using the theory of toric varieties. In fact, part of our aim is to show how the combinatorial result mentioned above can be naturally treated as a vanishing result for higher cohomology groups of certain line bundles on toric varieties associated to root systems (theorems C and D on Section \ref{resultssection}), which improves the previously known results for these varieties. This raises the question of new vanishing results for toric varieties in general (see ~\cite{qendrim} for work in this direction concerning toric varieties associated with root systems). It is worth mentioning that these toric varieties, where the fan of the variety is the Weyl fan, appear naturally in the work of group compactifications (see e.g. ~\cite{brion}, pp.187--206, and references therein). 

To give a flavor of the type of problem that we deal with here, it is very instructive to use a picture (see below) to illustrate the statement of Theorem E for n=3. Consider a toric variety $V=V_{SL_3}$ with initial lattice (i.e., lattice of co-characters) $L= \Z^3/\Z(1,1,1)$ and the fan consisting of six maximal cones as pictured below in dotted lines. The dual lattice (of characters) is $L^{\vee}= \{(x_1,x_2,x_3)\in \Z^3 | x_1+x_2+x_3=0\}.$

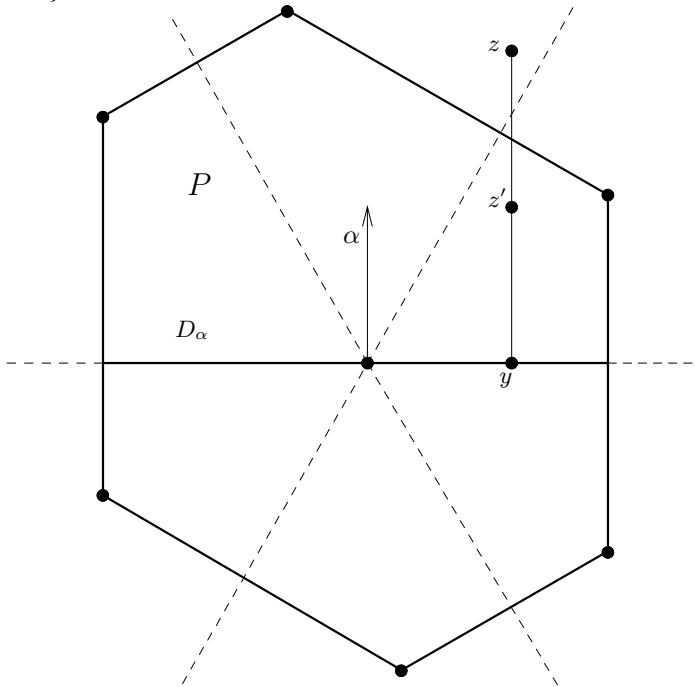
\setlength{\unitlength}{0.00083333in}
\begingroup\makeatletter\ifx\SetFigFont\undefined%
\gdef\SetFigFont#1#2#3#4#5{%
  \reset@font\fontsize{#1}{#2pt}%
  \fontfamily{#3}\fontseries{#4}\fontshape{#5}%
  \selectfont}%
\fi\endgroup%
{\renewcommand{\dashlinestretch}{30}
\begin{picture}(4299,4276)(0,-10)
\put(3162,2022){\blacken\ellipse{72}{72}}
\put(3162,2022){\ellipse{72}{72}}
\put(2262,2022){\blacken\ellipse{72}{72}}
\put(2262,2022){\ellipse{72}{72}}
\put(2472,102){\blacken\ellipse{72}{72}}
\put(2472,102){\ellipse{72}{72}}
\put(612,1197){\blacken\ellipse{72}{72}}
\put(612,1197){\ellipse{72}{72}}
\put(1762,4218){\blacken\ellipse{72}{72}}
\put(1762,4218){\ellipse{72}{72}}
\put(3762,3072){\blacken\ellipse{72}{72}}
\put(3762,3072){\ellipse{72}{72}}
\put(3762,840){\blacken\ellipse{72}{72}}
\put(3762,840){\ellipse{72}{72}}
\put(612,3558){\blacken\ellipse{72}{72}}
\put(612,3558){\ellipse{72}{72}}
\put(3162,3972){\blacken\ellipse{72}{72}}
\put(3162,3972){\ellipse{72}{72}}
\put(3162,2997){\blacken\ellipse{72}{72}}
\put(3162,2997){\ellipse{72}{72}}
\dashline{60.000}(1046,4169)(2262,2022)(3446,12)
\dashline{60.000}(3544,4242)(1107,20)
\thicklines
\path(3762,3072)(1764,4225)
\thinlines
\dashline{60.000}(3762,2022)(4287,2022)
\thicklines
\path(3769,852)(2470,102)
\path(612,1197)(2475,105)
\thinlines
\path(2292.000,2877.000)(2262.000,2997.000)(2232.000,2877.000)
\path(2262,2997)(2262,2021)
\path(3162,3972)(3162,2996)
\path(3162,2997)(3162,2021)
\thicklines
\path(604,3556)(1775,4231)
\path(612,3546)(612,1197)
\thinlines
\dashline{60.000}(612,2022)(12,2022)
\thicklines
\path(3762,3072)(3762,846)
\path(612,2022)(2262,2022)(3162,2022)(3762,2022)
\put(2112,2772){\makebox(0,0)[lb]{{\SetFigFont{10}{12.0}{\rmdefault}{\mddefault}{\updefault}$\alpha$}}}
\put(3087,1872){\makebox(0,0)[lb]{{\SetFigFont{9}{10.8}{\rmdefault}{\mddefault}{\updefault}$y$}}}
\put(1062,2172){\makebox(0,0)[lb]{{\SetFigFont{9}{10.8}{\rmdefault}{\mddefault}{\updefault}$D_{\alpha}$}}}
\put(1137,3072){\makebox(0,0)[lb]{{\SetFigFont{12}{14.4}{\rmdefault}{\mddefault}{\updefault}$P$}}}
\put(3012,2997){\makebox(0,0)[lb]{{\SetFigFont{9}{10.8}{\rmdefault}{\mddefault}{\updefault}$z'$}}}
\put(3012,3972){\makebox(0,0)[lb]{{\SetFigFont{9}{10.8}{\rmdefault}{\mddefault}{\updefault}$z$}}}
\end{picture}
}

Let $\mathcal{L}$ be a globally-generated torus-equivariant line bundle on $V$. Then $\mathcal{L}$ is completely determined by a set of characters (one for each maximal cone in the fan), such that if $\sigma$ and $\sigma '$ are two ``neighboring'' maximal cones, and $u(\sigma)$ and $u(\sigma ')$ are their corresponding characters, then $u(\sigma)-u(\sigma ')$, with respect to the canonical pairing of $L^{\vee}$ and $L$, is perpendicular to the common face of $\sigma$ and $\sigma '$, non-negative on $\sigma$, and non-positive on $\sigma '$.

This way $\mathcal{L}$ determines a hexagon $P$ where the vertices of $P$ are the characters corresponding to $\mathcal{L}$, one for each maximal cone of the fan. (In our picture all these characters are distinct and we have a true hexagon, but in general we may get some degenerate version of $P$, where some of the vertices coincide.) Then $\dim H^0(V, \mathcal{L})$ is given by the number of $L^{\vee}$-points in the convex hull $\Conv(P)$ of $P$. Note that $\Conv(P)$ lies in $L^{\vee}\otimes_{\Z}\R$. This should explain why we are using dotted lines for the cones of our initial fan --- they lie in the initial world, not the dual one.   

Fix a root $\alpha$ as in the picture. Then the ``line-segment" $D_{\alpha}$ corresponds to a (non-torus-invariant) divisor on $V$. If we denote by $i$ the embedding $D_{\alpha} \hookrightarrow V$, then $\dim H^0(V,i_*(\mathcal{L}|_{D_{\alpha}}))$ equals the number of points that are both on $D_{\alpha}$ and that are projections along the root $\alpha$ of points in $L^{\vee}$.

Our question is as follows: For any point $y$ on the line segment $D_{\alpha}$ that is obtained by the projection along $\alpha$ of some point $z \in L^{\vee}$ ($z$ is not necessarily in $\Conv(P)$), can we always find a point $z'$ that is in both $L^{\vee}$ and $\Conv(P)$, and which maps to the initial point $y$? The answer is ``yes" and, because of the long-exact sequence $$\cdots \longrightarrow H^0(V, \mathcal{L}) \stackrel{\varphi}{\longrightarrow} H^0(V,i_*(\mathcal{L}|_{D_{\alpha}})) \longrightarrow H^1(V,\mathcal{J}_{D_{\alpha}} \otimes \mathcal{L}) \longrightarrow 0,$$ where $\mathcal{J}_{D_{\alpha}}$ stands for the ideal sheaf of $D_{\alpha}$, it is equivalent to the statement that $H^1(V,\mathcal{J}_{D_{\alpha}} \otimes \mathcal{L})=0.$ We should mention here that the above long-exact sequence has trivial entries to the right of $H^1(V,\mathcal{J}_{D_{\alpha}} \otimes \mathcal{L})$, because higher cohomology groups vanish for globally-generated line bundles on complete toric varieties, and $\mathcal{L}$ and $\mathcal{L}|_{D_{\alpha}}$ are globally-generated on (the complete toric varieties) $V$ and $D_{\alpha}$, respectively. Also, note that the map $\varphi$ is induced by projecting along $\alpha$.

In this paper we precisely formulate the above question, using co-characters and co-roots, for a split connected reductive group $G$ over a finite extension of $\mathbb{Q}_p$, or equivalently, by passing to the (Langlands) complex dual world, using characters and roots, for the corresponding group $\hat{G}$ over $\mathbb{C}$. Then we reformulate the question using certain toric varieties associated to the root system of $\hat{G}$ (where the fan of our variety is the Weyl fan and the dual lattice (i.e., lattice of characters) is that of characters for $\hat{G}$). We answer the question in the affirmative for $G=GL_n$ and for a weaker version (where only line bundles $\mathcal{L}$ arising from Weyl orbits are considered) for $G=G_2$.

As mentioned at the beginning, this type of combinatorial question, without reference to toric varieties, was investigated by Kottwitz and Rapoport (and others). If we use the language of toric varieties to describe what they did, then they considered only globally generated line bundles $\mathcal{L}$ arising from Weyl-orbits, i.e., where the corresponding polygon $P$ is a Weyl-orbit for some element in the character lattice.

Let us now describe how our paper is organized. In the first section we introduce notation, set up the problem and state our main results (theorems A through E). We also show that Theorem C is equivalent to Theorem A. Sections 2 and 3 are devoted to some useful results in computing cohomology groups for line bundles on our toric varieties using the so-called piece-wise linear continuous functions on the support of the respective fans. We use these results in the next section, where we prove theorems B, D, and E, and where we show how E implies D and D implies B. Finally, the last section is devoted to the proof of Theorem C.

\paragraph{\bf{Acknowledgments:}} The author is greatly indebted to his thesis adviser, Professor Robert Kottwitz, for suggesting that he work on the problems dealt with in this paper, for numerous helpful comments, and for continuous encouragement and support.

\section{Set Up and Main Results}\label{resultssection}

Prior to ~\cite{cathy}, a converse to Mazur's inequality for $GL_{n}$ and $GSp_{2n}$ was proved in ~\cite{krapo}, Theorem 4.11 (see also ~\cite{frapo}, Theorem 2, for the $GL_{n}$ case), whose notation we follow. Once a notation is introduced, it will be fixed for the rest of the paper.

Let $F$ be a finite extension of $\mathbb{Q}_p$. Denote by $\mathcal{O}_F$ the ring of integers of $F$. Suppose $G$ is a split connected reductive group, $B$ a Borel subgroup and $T$ a maximal torus in $B$, all defined over $\mathcal{O}_F$. Let $P=MN$ be a parabolic subgroup of $G$ which contains $B$, where $M$ is the unique Levi subgroup of $P$ containing $T$.

We write $X$ for the set of co-characters $X_*(T)$. Then $X_{G}$ and $X_{M}$ will stand for the quotient of $X$ by the co-root lattice for $G$ and $M$, respectively. Also, we let $\varphi_{G}:X\rightarrow X_G$ and $\varphi_{M}:X\rightarrow X_M$ denote the respective natural projection maps.

Let $\mu\in X$ be $G$-dominant and let $W$ be the Weyl group of $T$ in $G$. The group $W$ acts on $X$ and so we consider $W\mu :=\left\{w(\mu): w\in W\right\}$ and the convex hull of $W\mu$ in $\mathfrak{a} :=X\otimes_{\mathbb{Z}}\mathbb{R}$, which we denote by $\Conv\left(W\mu\right)$. Define
$$P_{\mu}=\left\{\nu\in X: (i)\, \varphi_{G}(\nu)=\varphi_{G}(\mu); \, \emph{and} \,\,(ii)\, \nu\in \Conv\left(W\mu\right)\right\}.$$

Let $\mathfrak{a}_M:=X_M\otimes_{\Z}\R$ and write $\emph{pr}_M:\mathfrak{a}\rightarrow\mathfrak{a}_M$ for the natural projection induced by $\varphi_M$. Note that although $X_M$ is a quotient of $X$, after tensoring with $\mathbb{R}$ any possible torsion is lost and we can therefore consider $\mathfrak{a}_M$ as a subspace of $\mathfrak{a}$; we shall do so throughout the paper.

The aim is to generalize the proposition below for $G=GL_n$ and prove an analogous result for $G=G_2$. Kottwitz in ~\cite{kot} (sections 4.3 and 4.4) explains how, in each case, a converse to Mazur's inequality follows from the result below and therefore from our results.

\begin{proposition}\label{cathy}\emph{(cf. Theorem 0.2 in \cite{cathy})} Let $G$ be a split connected reductive group over $F$ with every irreducible component of its Dynkin diagram of type $A_n, B_n, C_n$, or $D_n$. With the notation as above, we have

$$\varphi_M\left(P_\mu\right)=\left\{\nu_1\in X_M :(i)\,\nu_1, \mu \,\,\emph{have the same image in}\,\, X_G;\right.$$
$$\left. \hspace{5cm} (ii)\, \emph{the image of}\, \,\nu_1\, \emph{in} \,\mathfrak{a}_M \,\emph{lies in}\, pr_M\left(\Conv W\mu \right)\right\}.$$
\end{proposition}

We are first going to phrase our results using Arthur's notion of $(G,T)$-orthogonal sets (see e.g. \cite{arthur}, pg. 217). For more on these sets see for example \cite{akot}, pp.441--447, whose notation we follow. A family of points $x_B$ in $X$, one for every Borel subgroup $B$ of $G$ that contains $T$, is called a \emph{$(G,T)$-orthogonal set} in $X$ if for every two adjacent Borel subgroups $B,B'$, there exists an integer $n$ such that $$x_B-x_{B'}=n {\alpha}^{\vee},$$ where $\alpha^{\vee}$ is the unique co-root for $T$ that is positive for $B$ and negative for $B'$. Similarly, if we have a point $x_P \in X_M$ for each parabolic subgroup $P$ of $G$ that admits $M$ as a Levi component, then the family $(x_P)$ is called a \emph{$(G,M)$-orthogonal set} in $X_M$ (see e.g. ~\cite{akot}, pg.442) if for every two adjacent parabolic subgroups $P, P'$ that admit $M$ as a Levi component, there exists an integer $m$ so that $$x_P-x_{P'}=m {\beta}^{\vee}_{P,P'}.$$ Here $\beta^{\vee}_{P,P'}$ is the unique element in $\Gamma$ with the property that all the other elements of $\Gamma$ are positive multiples of it; $\Gamma$ consists of the images in $X_M$ of the co-roots $\alpha^{\vee}$ where $\alpha$ is a root occurring in Lie$(N)\, \cap$ Lie$(\bar{N'})$, where $P=MN$ and $\bar{P'}=M\bar{N'}$, with $\bar{P'}$ being the parabolic subgroup containing $M$ opposite $P$.   

If all the numbers $n$ (resp. $m$) above are non-negative then we say that $(x_B)$ is a \emph{positive} $(G,T)$-orthogonal set (resp. $(x_P)$ is a \emph{positive} $(G,M)$-orthogonal set).

An important example of a positive $(G,T)$-orthogonal set arises from Weyl orbits $W\mu$ (see e.g. \cite{akot}, pg.443). Let $\nu \in X$, then we say that $\nu$ is dominant with respect to a Borel group $ B=TN $ if for every root $\alpha$ in Lie$(N)$ we have that $\left\langle \alpha,\nu  \right\rangle \geq 0.$ We get a positive $(G,T)$-orthogonal set by associating to every Borel group $B$ the unique element $x_B \in X$ that is both dominant with respect to $B$ and lies in $W\mu$.  

If $(x_B)$ is a $(G,T)$-orthogonal set in $X$, then we get a $(G,M)$-orthogonal set as follows. The points $(x_B)$, where $B$ is a Borel subgroup containing $T$ and $B \subset P$, form an $(M,T)$-orthogonal set in $X$, and we get a point $x_P\in X_M$ as the common image in $X_M$ of all the points in $\{x_B:B \subset P \}$. The set of all such points $x_P$, where $P$ is a parabolic that contains $B$ and admits $M$ as a Levi component, is a $(G,M)$-orthogonal set in $X_M$. Moreover, if $(x_B)$ is positive, then $(x_P)$ obtained in this way is positive as well.

We can now state our main results.
\newtheorem*{thma}{Theorem A}
\begin{thma} Let $G=G_2$. With notation as above we have that $$pr_M\left( \Conv(x_B)\cap X \right) = \Conv(x_P) \cap pr_M(X),$$
for every $(G,T)$-orthogonal set $(x_B)$ that arises from a Weyl orbit and its corresponding $(G,M)$-orthogonal set $(x_P)$.

\end{thma}

\newtheorem*{thmb}{Theorem B}
\begin{thmb} Let $G=GL_n$. With notation as above we have that $$pr_M\left( \Conv(x_B)\cap X \right) = \Conv(x_P) \cap pr_M(X),$$
for every positive $(G,T)$-orthogonal set $(x_B)$ and its corresponding $(G,M)$-orthogonal set $(x_P)$. 

\end{thmb}  

It is clear that the $GL_n$ case of Proposition \ref{cathy} becomes a special case of Theorem B when $(x_B)$ arises from a Weyl orbit. Also note that in both theorems A and B, the left-hand side is obviously contained in the right-hand side. The non-trivial part is to show the other containment.

It is important to mention that, while we believe that Theorem A is also true for the rest of the exceptional groups, Theorem B is probably only true for a split connected reductive group with irreducible components of Dynkin diagram of type $A_n$. More precisely, it is easy to construct counterexamples to Theorem B for other classical groups and for the $G_2$ case. At the end of section 5 we do so explicitly, but only for $G_2$, since the construction for classical groups is completely analogous. We believe that the same construction should also provide counterexamples to Theorem B for the other exceptional groups.

The proof of these two theorems will involve the theory of toric varieties and we will freely use standard terminology and basic facts from this theory which appear in ~\cite{fulton}. Let us first make clear the correspondence between $(G,T)$-orthogonal sets and line bundles on a certain projective nonsingular toric variety $V_G$, which we now define; we follow the exposition in ~\cite{akot}, \S 23, and keep the same notation as above. Let $\hat{G}$ and $\hat{T}$ be the (Langlands) complex dual group for $G$ and $T$, respectively. Let $Z(\hat{G})$ be the center of $\hat{G}$. Then the fan of our toric variety $V_G$ is the Weyl fan in $X_* (\hat{T}/Z(\hat{G})) \otimes_{\Z}\R$ and the torus is $\hat{T}/Z(\hat{G})$. (We would like to remark that the toric variety $V_G$ appears naturally in the theory of group compactifications---see e.g. ~\cite{brion}, pp.187--206.)

Clearly $\hat{T}/Z(\hat{G})$ acts on $V_G$, but we are interested in the action of $\hat{T}$ on $V_G$. The latter action is obtained using the canonical surjection $\hat{T} \twoheadrightarrow \hat{T}/Z(\hat{G})$ and the action of $\hat{T}/Z(\hat{G})$ on $V_G$.

As with any toric variety, there is a one-to-one correspondence between $\hat{T}$-orbits in $V_G$ and cones in the fan of $V_G$. In our case, because the fan is the Weyl fan, this means that we have a one-to-one correspondence between $\hat{T}$-orbits in $V_G$ and parabolic subgroups of $G$ that contain $T$. According to this identification, if $P$ is a parabolic subgroup of $G$ that contains $T$, then the corresponding $\hat{T}$-orbit is given by $\hat{T}/Z(\hat{M})$, where $M$ is the unique Levi component of $P$ containing $T$ and $\hat{M}$ is the corresponding Levi subgroup of $\hat{G}$ containing $\hat{T}$.

Again, as with any toric variety, there is a one-to-one correspondence between maximal cones in the fan of $V_G$ and $\hat{T}$-fixed points in $V_G$. In our case, since maximal cones correspond to Borel subgroups of $G$ that contain $T$, we get a one-to-one correspondence between Borel subgroups of $G$ containing $T$ and $\hat{T}$-fixed points in $V_G$.

It is very important for us to note that there is a one-to-one correspondence between isomorphism classes of $\hat{T}$-equivariant line bundles on $V_G$ and $(G,T)$-orthogonal sets in $X$ (see e.g. ~\cite{akot}, \S23.1). We describe the map that gives this correspondence. Let $\mathcal{L}$ be a $\hat{T}$-equivariant line bundle on $V_G$. Then at each $\hat{T}$-fixed point $v$ in $V_G$ the torus $\hat{T}$ acts by a character, say, $x_B$ on the line in $\mathcal{L}$ at $v$, where $B$ is the Borel subgroup corresponding to the fixed point $v$. In this way, for every Borel subgroup $B$ (containing $T$) we get a character $x_B\in X^*(\hat{T})$, i.e., a co-character $x_B \in X $. The fact that $(x_B)$ is a $(G,T)$-orthogonal set comes from the fact that the characters defining $\mathcal{L}$ must agree on the overlaps.

A very useful remark is that positive $(G,T)$-orthogonal sets correspond to line bundles on $V_G$ which are generated by their sections, and hence their higher cohomology groups vanish. This last result follows from the more general fact that if the support of the fan (i.e., the union of all the cones in the fan) of a toric variety is a convex set, then the higher cohomology groups vanish for line bundles (on this variety) which are generated by their sections (see e.g. ~\cite{fulton}, pg.74).

Remember that we are trying to reformulate theorems A and B in terms of toric varieties; we have already described the toric analog of $(G,T)$-orthogonal sets; the map $pr_M$ corresponds to the map $$X^*(\hat{T}) \twoheadrightarrow X^*(\hat{T})/R_{\hat{M}},$$ where $R_{\hat{M}}$ stands for the root lattice for $\hat{M}$ (note that the codomain of the last map is just $X^*(Z(\hat{M}))$); now we need to describe the toric analog of $(G,M)$-orthogonal sets. 

We continue to have the same notation. Then for $M$ as above (a Levi subgroup containing $T$), we will need a toric variety $Y_M^G$ for the torus $Z(\hat{M})/Z(\hat{G})$ (see e.g. ~\cite{akot}, \S23.2). First assume that the group $\hat{G}$ is adjoint; we can do this since the center of $\hat{M}/Z(\hat{G})$ in $\hat{G}/Z(\hat{G})$ is equal to $Z(\hat{M})/Z(\hat{G})$). Then $Z(\hat{M})$ is a subtorus of $\hat{T}$ and so $X_*(Z(\hat{M}))$ is a subgroup of $X_*(\hat{T})$. The collection of cones from the Weyl fan inside $X_*(\hat{T})\otimes_{\mathbb{Z}}\mathbb{R}$ that lie in the subspace $X_*(Z(\hat{M})) \otimes_{\mathbb{Z}}\mathbb{R}$ gives a fan. This is the fan for the complete, nonsingular, projective toric variety $Y_M^G$.

Now we go back to the general case where we no longer assume that $\hat{G}$ is adjoint. Denote by $G_0$ the simply connected cover of the derived group of $G$ and by $M_0$ the Levi subgroup in $G_0$ that is obtained as the inverse image of $M$ under the map $G_0 \rightarrow G$. Then we have that $\hat{G_0}=\hat{G}/Z(\hat{G})$, $\hat{M_0}=\hat{M}/Z(\hat{G})$, and $Z(\hat{M_0})=Z(\hat{M})/Z(\hat{G})$. From the adjoint case, we can define the toric variety $Y_{M_0}^{G_0}$ for $Z(\hat{M_0})=Z(\hat{M})/Z(\hat{G})$. We can therefore view $Y_{M_0}^{G_0}$  as a space on which $Z(\hat{M})$ acts using the canonical map $Z(\hat{M}) \twoheadrightarrow Z(\hat{M_0})$. It is this toric variety, which sits inside $V_G$ as a closed $Z(\hat{M})$-stable subspace, that we denote by $Y_M^G$. 

As with the toric variety $V_G$, it is easily seen that $Z(\hat{M})$-fixed points in $Y_M^G$ are in one-to-one correspondence with parabolic subgroups of $G$ that admit $M$ as a Levi component. Also, there is a one-to-one correspondence between isomorphism classes of $Z(\hat{M})$-equivariant line bundles on $Y_M^G$ and $(G,M)$-orthogonal sets in $X_M$. This correspondence is defined is the same way as the one for $V_G$ (see e.g. ~\cite{akot}, \S23.4).

We note that (see e.g. ~\cite{akot}, \S23.4 ) restricting a $\hat{T}$-equivariant line bundle $\mathcal{L}$ on $V_G$ to $Y_M^G$ corresponds the procedure described earlier in this section of obtaining a $(G,M)$-orthogonal set (the one corresponding to the restriction $\mathcal{L}|_{Y_M^G}$) from a $(G,T)$-orthogonal set (the one corresponding to $\mathcal{L}$).

Let us now assume that our parabolic subgroup $P$ is of semisimple rank 1. This implies that the root lattice $R_{\hat{M}}$ is just $\mathbb{Z} \alpha$ for a unique, up to a sign, root $\alpha$ of $\hat{G}$, and that the toric variety $Y_M^G$, which we now denote by $D_{\alpha}$, is a non-torus-invariant divisor in $V_G$. The map $pr_M$ will now be denoted by $p_{\alpha}$. By tensoring with $\mathbb{R}$ we get a map from $p_{\alpha}$, which we still denote by $$p_{\alpha}:X^*(\hat{T}) \otimes_{\mathbb{Z}}\mathbb{R} \twoheadrightarrow ( X^*(\hat{T})/\mathbb{Z}\alpha )\otimes_{\mathbb{Z}}\mathbb{R} .$$ Since tensoring with $\mathbb{R}$ will loose any possible torsion, we can identify the codomain of the last map with the co-root hyperplane $$[\alpha^{\vee} = 0] : =\{x \in  X^*(\hat{T})\otimes_{\mathbb{Z}}\mathbb{R}: \langle \alpha^{\vee},x \rangle =0  \} ,$$ where $\langle ,\rangle$ is the canonical pairing between co-characters and characters, and $\alpha^{\vee}$ stands for the co-root of $\hat{G}$ corresponding to $\alpha$. We will therefore say that the map $p_{\alpha}$ is the projection along $\alpha$ onto the co-root hyperplane $[\alpha^{\vee} = 0]$. Also, the fan of $D_{\alpha}$ is contained in the root hyperplane $$[\alpha = 0] : = \{ x \in X_*(\hat{T}) \otimes_{\mathbb{Z}}\mathbb{R}: \langle x, \alpha \rangle =0 \}.$$

Now let $\mathcal{L}$ be a $\hat{T}$- line bundle on $V_G$ that is generated by its sections. Then we have a short-exact sequence of sheaves on $V_G$: $$0\longrightarrow \mathcal{J}_{D_{\alpha}} \otimes \mathcal{L} \longrightarrow \mathcal{L} \longrightarrow i_*(\mathcal{L}|_{D_{\alpha}}) \longrightarrow 0,$$
where $\mathcal{J}_{D_{\alpha}}$ is the ideal sheaf of $D_{\alpha}$ and $i$ is the inclusion map $D_{\alpha} \hookrightarrow V_G$. Note that $H^i(V_G, \mathcal{L})$=0, for all $i>0$, since $\mathcal{L}$ is generated by its sections and $V_G$ is projective, and also $$H^i(V_G,i_*(\mathcal{L}|_{D_{\alpha}}))=H^i(D_{\alpha},\mathcal{L}|_{D_{\alpha}})=0,$$ for all $i>0$, since $\mathcal{L}|_{D_{\alpha}}$ is generated by its sections and $D_{\alpha}$ is a projective toric variety. Therefore the short-exact sequence gives rise to the long-exact sequence
$$(*)\,\,\,\,\,... \longrightarrow H^0(V_G, \mathcal{L}) \stackrel{\varphi}{\longrightarrow} H^0(V_G,i_*(\mathcal{L}|_{D_{\alpha}})) \longrightarrow H^1(V_G,\mathcal{J}_{D_{\alpha}} \otimes \mathcal{L}) \longrightarrow 0.$$ So, we see that surjectivity of the map $\varphi$ is equivalent to $H^1(V_G,\mathcal{J}_{D_{\alpha}} \otimes \mathcal{L})=0.$

Let $G=G_2$. We claim that the following result is equivalent to  Theorem A.

\newtheorem*{thmc}{Theorem C}
\begin{thmc}
With notation as above, we have that $$H^1(V_{G_2},\mathcal{J}_{D_{\alpha}} \otimes \mathcal{L})=0,$$ whenever $\mathcal{L}$ arises from a Weyl orbit. 
\end{thmc}

\begin{remark}
Theorem C is true for split classical groups. In fact, in the same way that we will show the equivalence of theorems A and C one can show that theorem C for a split classical group is equivalent to Proposition \ref{cathy} for the case of parabolic subgroups $P$ of semisimple rank 1.
\end{remark}

Let us demonstrate the equivalence between theorems A and C. Put $G=G_2$. Suppose that $\mathcal{L}$ corresponds to the Weyl-orbit orthogonal set $(x_B)$ in $X$, i.e. in $X^*(\hat{T})$, and, as before, denote by $(x_P)$ the corresponding $(G,M)$-orthogonal set in $X_M$, i.e. in $X^*(Z(\hat{M}))$. Write $P_{\mathcal{L}}$ for the intersection of $X^*(\hat{T})$ with the convex hull $\Conv\{x_B\}$ of $(x_B)$. Then we have (see e.g. ~\cite{fulton}, pg.66) $$H^0(V_G,\mathcal{L})= \bigoplus_{u\in P_{\mathcal{L}}} \C \, \chi^u,$$ where the section $ \chi^u$ is an eigenvector corresponding to the character $u$. (Recall that $\hat{T}$ acts on $H^0(V_G,\mathcal{L})$ and so $H^0(V_G,\mathcal{L})$ decomposes according to the characters of $\hat{T}$.) 

Assume that $\alpha$ is the root of $\hat{G}$, up to a sign, that corresponds to $\hat{M}$, in the sense that the projection map $p_{\alpha}: X^*(\hat{T}) \rightarrow [\alpha^{\vee} = 0]$ along the root $\alpha$ corresponds to the projection map $pr_M$. Then we write $P_{\mathcal{L},\alpha}$ for the intersection of $\Conv\{x_P\}$ with $p_{\alpha}(X^*(\hat{T}))$ and we get $$H^0(V_G,i_*(\mathcal{L}|_{D_{\alpha}}))= \bigoplus_{u \in P_{\mathcal{L},\alpha}}\C \, \chi^u.$$
The map $\varphi:H^0(V_G,\mathcal{L}) \rightarrow H^0(V_G,i_*(\mathcal{L}|_{D_{\alpha}}))$ from the long-exact sequence above is given by $\varphi(\chi^u)=\chi^{p_{\alpha}(u)}$. This explains why Theorem C is equivalent to Theorem A, since $p_{\alpha}$ corresponds to $pr_M$.

Now let $G=GL_n$. Theorem B will follow from the following result.

\newtheorem*{thmd}{Theorem D}
\begin{thmd}
With notation as above, we have that $$H^1(V_{GL_n},\mathcal{J}_{D_{\alpha}} \otimes \mathcal{L})=0,$$ whenever $\mathcal{L}$ is generated by its sections.

\end{thmd}

Theorem D itself will be deduced (in Section 4.5) from the result for $G=SL_n$:

\newtheorem*{thme}{Theorem E}
\begin{thme}
With notation as above, we have that $$H^1(V_{SL_n},\mathcal{J}_{D_{\alpha}} \otimes \mathcal{L})=0,$$ whenever $\mathcal{L}$ is generated by its sections.

\end{thme}

We note that while theorems A and C were seen to be equivalent, to show that theorems B and D imply each other will take a little more work (because we are no longer working with a two-dimensional space!); however, the same argument above for $G=G_2$ applies in our case to show that Theorem B for the case of parabolic subgroups $P$ of semisimple rank 1 is equivalent to Theorem D. But, we will deduce our Theorem B \emph{from} Theorem D, where the latter will be proved using toric geometry.

We explain in Section \ref{Proof of Theorem B} how the general case is handled, but here let us just show, using a similar method used there, how Theorem D, i.e. surjectivity of $\varphi$, implies Theorem B for all $M=GL_{n_1}\times GL_{n_2}\times \cdots \times GL_{n_r}$, where $\sum_{i=1}^{r}n_i = n$ and only one of $n_i$'s is equal to 2 and the rest are equal to 1.

Assume that $n_k=2$ for some $1 \leq k \leq n_r$ and that $n_i=1$ for all $i$ such that $i \neq k$ and $1 \leq i \leq n_r$. Keeping the same notation as above, we have in this case $\mathfrak{a}=\R^n$ and $$\mathfrak{a}_M=\{(x_1,\ldots,x_n)\in \R^n:x_k=x_{k+1}\},$$ where, as mentioned earlier, we have identified $\mathfrak{a}_M$ as a subspace of $\mathfrak{a}$. Also, the map $pr_M$, which in general is given by averaging over each of the $r$ ``batches'', now agrees with $p_{\alpha}$ where $\alpha=(0,\ldots,0,1,-1,0,\ldots,0)$ (with 1 in the $k$-th place and -1 in the $(k+1)$-st place): $$p_{\alpha}(x_1,\ldots,x_n)=(x_1,\ldots,x_{k-1},\frac{x_k+x_{k+1}}{2},\frac{x_k+x_{k+1}}{2},x_{k+2},\ldots,x_n).$$ (In general $pr_M$ does not correspond to one $p_{\alpha}$, rather, it is a composition of a number of $p_{\alpha}$'s, for distinct roots $\alpha$.)

As in the $G_2$ case, we find that $$H^0(V_G,\mathcal{L})= \bigoplus_{u\in P_{\mathcal{L}}} \C \, \chi^u,$$
$$H^0(V_G,i_*(\mathcal{L}|_{D_{\alpha}}))= \bigoplus_{u \in P_{\mathcal{L},\alpha}}\C \, \chi^u,$$
and $\varphi (\chi^u)=\chi^{p_{\alpha}(u)}$. Then surjectivity of $\varphi$ clearly implies Theorem B for this particular case.

As we mentioned in the Introduction, the next sections are devoted to proving our results. It is now safe for the reader to assume that, for the rest of the paper, we are working over the complex numbers---we will be in the toric geometry setting. Also, since we will be working in the complex dual world, when we use terms root, co-root, character, and co-character, they will always refer to those terms for the complex dual groups $\hat{G}$ and $\hat{T}$.

Before we end this section it is worth mentioning that vanishing results like the ones in theorems C and D (and E) are also of independent interest just from a toric-variety point of view. A very important vanishing result for toric varieties has been proved by Musta\c{t}\u{a} (see e.g. ~\cite{mustata}, Theorem 0.1) and for the particular toric varieties arising from $GL_n$ and $G_2$ theorems C and D (and E) give vanishings of higher cohomology groups for more line bundles on these varieties. (See also ~\cite{qendrim} for more vanishing results for toric varieties associated with root systems.)

\section{Some useful results}

We keep the same notation as in the previous section. Let $V=V_G$ be the non-singular, projective toric variety corresponding to $G$ as described in Section \ref{resultssection} and let $\Delta$ be the fan of $V$. Each element of $\Delta$ corresponds to a unique torus orbit in $V$, and under this identification, if we denote by $\Delta(1)$ the set of rays in $\Delta$ and if $\rho \in \Delta(1)$, we write $D_{\rho}$ for the closure of the orbit of $\rho$, a $\hat{T}$-Cartier divisor (cf. ~\cite{fulton},\S3.1 and \S3.3). In fact, $\hat{T}$-Cartier divisors on $V$ can be written as a linear combination of $D_{\rho}$'s. From now on let us agree to refer to $\hat{T}$-Cartier divisors simply as divisors. (No confusion should arise from this agreement even though our divisor $D_{\alpha}$ is not a $\hat{T}$-divisor; in fact, thanks to Lemma \ref{Dalpha2} below, when computing cohomology groups, we will deal with a $\hat{T}$-Cartier divisor instead of $D_{\alpha}$.)

Let us mention another way of characterizing divisors on toric varieties (see e.g. ~\cite{fulton}, \S3.4, but be aware of a sign difference). A divisor $D$ on a toric variety $V$, with fan $\Delta$, is completely determined by a set of characters $u(\sigma)$, one for every element $\sigma \in \Max(\Delta)$, such that they ``agree'' on the overlaps, i.e., such that they form a $(G,T)$-orthogonal set (we will usually call these just orthogonal sets; here and throughout this paper, $\Max(\Delta)$ stands for the set of maximal cones in $\Delta$). This orthogonal set in turn gives a continuous (integral) piece-wise linear function $\psi_{D}$ on the \emph{support} $\left| \Delta \right|:= \cup_{\sigma \in |\Delta|}\sigma$ of the fan, where $\psi_D(v)=\left\langle u(\sigma), v\right\rangle$, for all $v\in \sigma$. And conversely, given a continuous (integral) piece-wise linear function $\psi$ on $\left|\Delta \right|$, i.e., a continuous function that is linear and given by an element of the character lattice on each cone of $\Delta$, we get a divisor by considering the characters $u(\sigma)$, $\sigma \in \Max(\Delta)$, which define $\psi$ on the maximal affine pieces of $V$. (The continuity of $\psi$ guarantees that the characters defining it agree on the overlaps and hence give a divisor.)  

From what we said, it is easy to see that if $D=\sum_{\rho \in \Delta} a_{\rho}D_{\rho}$, then the corresponding piece-wise linear function is given by $\psi_{D}(v)=\left\langle u(\sigma), v\right\rangle$ $(v\in \sigma)$ where, for each $\sigma$, the character $u(\sigma)$ is found by solving the system of equations $$\left\langle u(\sigma), v_{\rho} \right\rangle = a_{\rho}\, ,$$ where $\rho$ varies through the rays in the cone $\sigma$ and $v_{\rho}$ stands for the primitive lattice element in the ray $\rho$.   

Let us mention that $\psi_D$ is called \emph{convex} or \emph{lower convex} if for any $v,v'\in |\Delta|$ and $0\leq p \leq 1$, we have that $\psi_D(pv+(1-p)v') \leq p\psi_D(v) + (1-p)\psi_D(v')$.

The following result is going to be very useful.

\begin{proposition}\label{fultontheorem}
\begin{itemize}
\item[a)] \emph{(cf. ~\cite{fulton}, pg. 68)} Let $D$ be a Cartier divisor on $V$. Then $\mathcal{O} (D)$ is generated by its sections if and only if $\psi_D$ is (lower) convex.

\item[b)] \emph{(cf. ~\cite{fulton}, pg.74)} For every $i \geq 0$ there exist canonical isomorphisms $$H^i(V, \mathcal{O}(D)) \cong \bigoplus_{u\in X^*(\hat{T})} H^i(V, \mathcal{O}(D))_u,$$
and $$H^i(V, \mathcal{O}(D))_u \cong H^i(\left|\Delta \right|,\left|\Delta \right| \setminus Z(u); \C),$$
where $H^i(V, \mathcal{O}(D))_u$ denotes the $u$-eigenspace that is obtained by the action of $\hat{T}$ on $H^i(V, \mathcal{O}(D))$ and $Z(u)=\{v\in |\Delta|: \psi(v)\geq \left\langle u, v \right\rangle \}$.

\end{itemize}
\end{proposition}

Keep the same notation as above and let $\alpha$ be a root of $\hat{G}$. The following result was suggested to us by R. Kottwitz.

\begin{lemma}\label{Dalpha} For the toric variety $V=V_G$, we have that $D_{\alpha}$ and $$\sum_{\rho \in \Delta(1),\,\left\langle \alpha, v_\rho \right\rangle >0 }\left\langle \alpha, v_\rho \right\rangle D_{\rho}$$ are linearly equivalent divisors.

\end{lemma}
\begin{proof}
This follows from considering the rational function $\chi^{\alpha} - 1$ on $V_G$ and computing the corresponding principal divisor. Note that the divisor of zeros corresponds to $D_{\alpha}$.
\end{proof}

Since we will be interested in computing the cohomology group $H^1(V_G,\mathcal{O}({D-D_{\alpha}}))$, the last lemma allows us  to focus instead on $$H^1(V_G,\mathcal{O}({D-\sum_{\rho \in \Delta(1),\,\left\langle \alpha, v_\rho \right\rangle >0 }\left\langle \alpha, v_\rho \right\rangle D_{\rho}})).$$

The lemma below describes the orthogonal set corresponding to the $\hat{T}$-equivariant divisor $ \sum_{\rho \in \Delta(1),\,\left\langle \alpha, v_\rho \right\rangle >0 }\left\langle \alpha, v_\rho \right\rangle D_{\rho} $, which, due to the previous lemma, will be referred to as the orthogonal set corresponding to $D_{\alpha}$.

\begin{lemma}\label{Dalpha2}
 Let $\alpha$ be a root of $\hat{G}$ and consider $V_G$. The corresponding orthogonal set for $D_{\alpha}$ is given by $\left( u(\sigma) \right)_{\sigma \in \Max(\Delta)}$, where $u(\sigma)=\alpha$ if $\alpha\geq0$ on $\sigma$, and $u(\sigma)=0$ if $\alpha \leq 0$ on $\sigma$.
\end{lemma}

\begin{proof} The result follows at once from Lemma \ref{Dalpha} and the fact that, as we discussed earlier, if $D=\sum_{\rho \in \Delta} a_{\rho}D_{\rho}$, then $u(\sigma)$'s can be found by solving the system of equations $\left\langle u(\sigma), v_{\rho} \right\rangle = a_{\rho}$, where $\rho \in \Delta(1)$ varies through the rays of $\sigma$.
\end{proof}

\section{An Important Lemma}\label{lemmasection}

\begin{lemma}\label{importantlemma} Consider $V_G$ and let $\alpha$ be a root of $\hat{G}$. Suppose $\mathcal{O}(D)$ is generated by its sections. If one of the conditions \emph{(i)} or \emph{(ii)} below is satisfied, then the 0-th eigenspace $$H^1(V_G,\mathcal{J}_{D_{\alpha}} \otimes \mathcal{O}(D))_0=0.$$

\begin{itemize}
\item [(i)]
$\left\{ v\in |\Delta|: \psi_D(v) < 0\, \emph{and} \left\langle \alpha, v \right\rangle \leq 0 \right\} = \varnothing;$

\item [(ii)] $\left\{ v\in |\Delta|: \psi_D(v) < 0\, \emph{and} \left\langle \alpha, v \right\rangle \geq 0 \right\} \neq \varnothing$ and 

$\left\{ v\in |\Delta|: \psi_D(v) < 0\, \emph{and} \left\langle \alpha, v \right\rangle \leq 0 \right\} \neq \varnothing.$

\end{itemize}
\end{lemma}
\begin{proof}
First note that using Proposition \ref{fultontheorem}, part b), it is sufficient (and necessary) to show that
\begin{equation}\label{eqn:vanish}
H^1(\left|\Delta \right|,U_{D,\alpha}(0);\C)=0,
\end{equation}
where $U_{D,\alpha}(0)=\left\{ v\in |\Delta|: \psi_{D-D_{\alpha}}(v) < 0 \right\}.$

Looking at the long-exact sequence we get from the pair $(\left|\Delta \right|,U_{D,\alpha}(0))$, since $|\Delta|$ is a vector space ($V_G$ is complete), we only need to show that $U_{D,\alpha}(0)$ is path-wise connected. Let us do that.

According to Lemma \ref{Dalpha2}, $U_{D,\alpha}(0)$ is equal to the union of sets
\begin{equation}\label{eqn:psiD}
U_{D,\alpha}(0)\cap [\alpha \leq 0] = \left\{ v\in |\Delta|: \psi_{D}(v) < 0\, \emph{and} \left\langle \alpha, v\right\rangle \leq 0 \right\}
\end{equation}
and 
\begin{equation}\label{eqn:psiDalpha}
U_{D,\alpha}(0)\cap [\alpha \geq 0] =\left\{ v\in |\Delta|: \psi_{D-D_{\alpha}}(v) < 0\, \emph{and} \left\langle \alpha, v\right\rangle \geq 0 \right\}.
\end{equation}
Here we have written $[\alpha \geq 0]$ and $[\alpha \leq 0]$ for the sets $\{v\in |\Delta| : \left\langle \alpha, v \right\rangle \geq 0 \}$ and $\{v\in |\Delta| : \left\langle \alpha, v \right\rangle \leq 0 \}$, respectively.

Proposition \ref{fultontheorem}, part a), implies that sets (\ref{eqn:psiD}) and (\ref{eqn:psiDalpha}) are convex. Let us see how this works in the case of (\ref{eqn:psiDalpha}). Let $v$ and $v'$ be in (\ref{eqn:psiDalpha}) and let $0 \leq p \leq 1$. Then using Lemma \ref{Dalpha2}. we have $$\psi_{D-D_{\alpha}}(pv+(1-p)v')= \psi_D(pv+(1-p)v')- \langle \alpha, pv+(1-p)v' \rangle,$$
which, since $\psi_D$ is convex, is not greater than
$$p\psi_D(v)+(1-p)\psi_D(v')-p \langle \alpha, v \rangle - (1-p) \langle \alpha, v' \rangle = p \psi_{D-D_{\alpha}}(v) + (1-p) \psi_{D-D_{\alpha}}(v'),$$
and the last expression is less than zero since $\psi_{D-D_{\alpha}}(v) < 0$, $\psi_{D-D_{\alpha}}(v') < 0$ and $0 \leq p \leq 1$. Clearly, from our assumptions, we also have that $pv+(1-p)v' \in [\alpha \geq 0]$, so (\ref{eqn:psiDalpha}) is a convex set, as desired.

If condition (i) of the lemma is satisfied, then $U_{D,\alpha}(0)$ equals the set (\ref{eqn:psiDalpha}) and hence (\ref{eqn:vanish}) holds, because the set (\ref{eqn:psiDalpha}) is convex.

Suppose that (ii) is satisfied. Since $\{v\in |\Delta|: \psi_D(v)<0 \}$ is convex, this implies that there exists a (non-empty, convex) subset $U$ in $ [\alpha \leq 0] \cap [\alpha \geq 0] $ that is contained in both sets appearing in (ii). We claim that $U$ is then contained in both sets (\ref{eqn:psiD}) and (\ref{eqn:psiDalpha}). This claim is sufficient to prove that $U_{D,\alpha}(0)$ is path-wise connected (and hence (\ref{eqn:vanish}) holds) since we already saw that (\ref{eqn:psiD}) and (\ref{eqn:psiDalpha}) are convex.

First, $U$ is contained in (\ref{eqn:psiD}) by our assumption that it is contained in both sets appearing in (ii), and in particular the second one. Second, since $U$ is contained in the first set appearing in (ii), we see that $\psi_D(v) < 0$ and $\langle \alpha,v \rangle \geq 0$, for all $v\in U$. Therefore, using Lemma \ref{Dalpha2}, we see that $$\psi_{D-D_{\alpha}}(v) = \psi_D(v) - \langle \alpha,v \rangle < 0.$$ Thus $U$ is contained in (\ref{eqn:psiDalpha}) as well, which ends the proof of the claim and hence of the lemma. 
\end{proof}

Note that if we knew that (\ref{eqn:vanish}) is true for \emph{all} globally generated divisors $D$, then we could conclude that for those divisors we have $H^1(V_G, \mathcal{O}(D-D_{\alpha}))_u = 0,$ for \emph{all} $u$. Indeed, to prove this, consider the (globally generated) divisor $D'$ associated to the orthogonal set $\{u(\sigma)-u\}$, where $\{u(\sigma)\}$ corresponds to $D$. The result follows by using (\ref{eqn:vanish}) (which is assumed to be true for all globally generated divisors), where instead of $U_{D,\alpha}(0)$ we have $U_{D',\alpha}(0)$.

Using Proposition \ref{fultontheorem}, part b), we deduce that if (\ref{eqn:vanish}) holds for all globally generated divisors $D$, then $H^1(V_G,\mathcal{J}_{D_{\alpha}} \otimes \mathcal{O}(D))=0$. 

An important consequence of this is that Theorem D needs to be proved only in the case when the set (\ref{eqn:psiD}) is non-empty and
\begin{equation}\label{eqn:psiDelta4}
\{ v\in |\Delta|: \psi_{D}(v) < 0\} \cap [\alpha \geq 0] = \varnothing.
\end{equation}
More concretely, it suffices to show that, under these conditions,
\begin{equation}\label{eqn:psiDelta5}
U_{D,\alpha}(0)\cap [\alpha \geq 0] = \varnothing ,
\end{equation}
 because then we get that $U_{D,\alpha}(0)$ is equal to the set (\ref{eqn:psiD}), which we know is a convex set, and we can apply the same reasoning as in Lemma \ref{importantlemma} to conclude that (\ref{eqn:vanish}) is true, and by what we just wrote, $H^1(V_G, \mathcal{J}_{D_{\alpha}} \otimes\mathcal{O}(D)=0$ is true as well.

\section{Proofs of Theorems B, D, and E}

Consider $V_G$, where $G={SL_n}$ (and so $\hat{G}=PGL(n,\mathbb{C})$). We want to concretely describe the fan $\Delta$ of $V_G$. But first, the initial lattice (i.e. lattice of characters) for $V_G$ is $L= \Z^n/\Z(1,1,\ldots,1)$ and so the dual one is $L^{\vee}= \{(x_1,\ldots,x_n)\in \Z^n | x_1+\ldots+x_n=0\}.$ For every $i=1,\ldots,n$ we put $L_i$ for (the element represented by) $(0,\ldots,0,1,0,\ldots,0)$, where 1 appears only on the $i$-th place. Identify the rays in $\Delta$ with their corresponding minimal lattice points. Then $\Delta (1)$ consists of the sums of the form $\sum^{k}_{j=1}L_{i_j}$, where $i_j$'s are distinct elements from $\{1,2,\ldots, n\}$ and $k=1,...,n-1$. The maximal cones in $V_G$ are the n-dimensional cones whose rays are of the type $$L_{i_1}, L_{i_1}+ L_{i_2},\ldots , L_{i_1}+L_{i_2}+\ldots +L_{i_{n-1}},$$ where $i_1, i_2,\ldots ,i_{n-1}$ are distinct elements from $\{1,2,\ldots, n \}$. Let $\alpha$ be a root of $\hat{G}$. We can assume that $\alpha $ is simple (for a choice of a Borel subgroup) and, with a change of coordinates, we can always ensure that $\alpha$ is written as $L_1-L_2$. Therefore, without loss of generality, throughout this section we let $\alpha = L_1-L_2$ .

The proof of Theorem E will be proved by induction, although we begin by giving the proof for the first two nontrivial cases, $n=3$ and $n=4$. This may seem a little unusual, but we do so because in order to apply the induction process we will need ``many" root hyperplanes and for $n=4$ we don't have ``enough" of them (see Case 2) of the proof in Section 4.3 where we use $n>4$).

\subsection{Proof of Theorem E for $n=3$}

In this case the set of maximal cones $\Max(\Delta)$ in the fan $\Delta$ of $V_G$ consists of six cones $\sigma_1,\ldots,\sigma_6$, generated respectively by $L_1$ and $-L_3$; $-L_3$ and $L_2$; $L_2$ and $-L_1$; $-L_1$ and $L_3$; $L_3$ and  $-L_2$; and $-L_2$ and $L_1$.

Let $D$ be a divisor on $V_G$ such that $\mathcal{O}(D)$ is generated by its sections. If $u(\sigma_i), i=1,\ldots,6$ form the orthogonal set corresponding to $D$, where $u(\sigma_i)$ is assigned to the maximal cone $\sigma_i$, then, according to Lemma \ref{Dalpha2}, the orthogonal set corresponding to the divisor $D-D_{\alpha}$ consists of $u(\sigma_i) - \alpha$ for $i=1,5,6$ and $u(\sigma_i)$ for $i=2,3,4$.   

From our discussion in Section \ref{lemmasection}, we only need to prove that (\ref{eqn:psiDelta5}) is true, under the assumption that the set (\ref{eqn:psiD}) is non-empty and (\ref{eqn:psiDelta4}) is true. Suppose, for a contradiction, that under these assumptions there exists $v_0 \in |\Delta|$ such that $\left\langle \alpha, v_0 \right\rangle \geq 0$ and $\psi_D(v_0)-\left\langle \alpha, v_0 \right\rangle < 0$. Since $\psi_D$ is continuous and piecewise linear, there are only two cases we need to consider: (i) $v_0=L_1$, and (ii) $v_0=L_1+L_3=-L_2$. We can use the symmetry $(x,y,z)\mapsto (-y,-x, -z)$ of the root system, which preserves the half-space  $[\alpha \geq 0]$, to see that we only need to consider the case (i).

Let $v_0 = L_1$. Suppose that $u(\sigma_1)=(a,b,-a-b)$, where $a,b\in \Z$. Since $\psi_D \geq 0$ on $[\alpha \geq 0]$, we find that $\left\langle u(\sigma_1), L_1\right\rangle \geq 0$ and $\left\langle u(\sigma_1), -L_3\right\rangle \geq 0$. Therefore $a \geq 0$ and $a+b\geq 0$. But, by assumption $\left\langle u(\sigma_1)-\alpha, L_1\right\rangle < 0$, hence $a-1 <0$. So $u(\sigma_1)= (0,b,-b)$ and $b\geq 0$. Now, since $\mathcal{O}(D)$ is generated by its sections, $u(\sigma_i)$'s form a \emph{positive} orthogonal set and therefore there exists a non-negative number $n\in \Z$ such that $u(\sigma_2)=u(\sigma_1)+n(L_2-L_1)$. We get $\psi_D(L_2)=\left\langle u(\sigma_2), L_2 \right\rangle =b+n \geq 0$. For the same reasons, there exists a non-negative integer $m$ such that $u(\sigma_3)=u(\sigma_2)+m(L_3-L_1)$, and therefore $\psi_D(-L_1)=\left\langle u(\sigma_3),-L_1 \right\rangle = n+m \geq 0$. Thus, since $\psi_D$ is continuous and piece-wise linear, $\psi_D$ takes non-negative values on all of $[\alpha \leq 0]$, and this contradicts the non-emptiness of (\ref{eqn:psiD}). Theorem E follows for $n=3$.

\begin{remark}\label{inductionremarkn=3}
If we are working with $V_{SL_4}$, then we can consider sub-fans of the initial fan $\Delta$ that are contained in the root hyperplanes. For example, if our given root is $\beta = L_3-L_4$, then the sub-fan corresponding to the root hyperplane $[\beta =0]$ has rays $L_1, -L_3-L_4, L_2, -L_1, L_3+L_4, -L_2$ and the maximal-dimensional cones are the 2-dimensional cones in $\Delta$ obtained from these rays. This is just a ``copy'' of the fan of $V_{SL_3}$ and indeed we can carry out the same calculations as above to see that Theorem E holds for this ``copy'' of $V_{SL_3}$ inside $V_{SL_4}$. This remark is important for our induction process and is part of a more general story which we tell in Section \ref{inductionsection}.   
\end{remark}

\subsection{Proof of Theorem E for $n=4$}

Note that we are still assuming that $\alpha = L_1-L_2$. To simplify notation, write $\sigma_{i,j,k}$ for the maximal cone whose rays are $L_i, L_i+L_j,L_i+L_j+L_k$, where $i,j,k$ are distinct elements of $\{ 1,2,3,4 \}$.

Let $D$ be a divisor on $V_G$ with $\mathcal{O}(D)$ generated by its sections. Denote by $u(\sigma_{i,j,k})$ the character corresponding to $D$ for the cone $\sigma_{i,j,k}$. Then these characters form an orthogonal set and, according to Lemma \ref{Dalpha2}, the orthogonal set corresponding to $D-D_{\alpha}$ consists of $u(\sigma_{i,j,k}) - \alpha $ for the maximal cones $\sigma_{i,j,k}$ contained in the half-space $[\alpha \geq 0]$ and $u(\sigma_{i,j,k})$ for the ones contained in $[\alpha \leq 0]$. Following the discussion in Section \ref{lemmasection}, we assume that the set (\ref{eqn:psiD}) is non-empty and that (\ref{eqn:psiDelta4}) is true. We want to prove (\ref{eqn:psiDelta5}) and, for a contradiction, suppose that there exists a point $v_0$ in the half-space $[\alpha \geq 0]$ such that $\psi_D(v_0)- \left\langle \alpha, v_0 \right\rangle <0 $. Since $\psi_D$ is continuous and piece-wise linear on $|\Delta|$, there are only four cases we need to consider, corresponding to the rays whose primitive lattice points are contained in $[\alpha >0]$:
(i) $v_0=L_1$, (ii) $v_0=L_1+L_3$, (iii) $v_0=L_1+L_4$, and (iv) $v_0=L_1+L_3+L_4$. But, we can use the symmetry $(x,y,z,w) \mapsto (-y,-x,-z,-w)$ of the root system, which preserves the half-space $[\alpha \geq 0]$, to see that we only need to consider cases (i) and (ii).

\emph{Case} (i): Our aim is to show that $\psi_D$ takes non-negative values at $L_2, L_2+L_3,L_2+L_4$ and $L_2+L_3+L_4$, which would contradict our assumption that the set (\ref{eqn:psiD}) is non-empty. First, if we look at the hyperplane corresponding to the root $L_2-L_3$, as mentioned in Remark \ref{inductionremarkn=3}, we get a copy of $V_{SL_3}$, where the rays of the sub-fan are $L_1,L_1+L_4, L_4, L_2+L_3+L_4,L_2+L_3, L_1+L_2+L_3$. These are precisely the rays of $\Delta$ contained in our hyperplane. Apply Theorem E, case n=3, to see that $\psi_D(L_2+L_3) \geq 0$ and $\psi_D(L_2+L_3+L_4) \geq 0$.

Similarly, by looking at the hyperplanes corresponding to $L_3-L_4$ and $L_2-L_4$ and applying Theorem E, for n=3, we get that $\psi_D(L_2) \geq 0$ and $\psi_D(L_2+L_4) \geq 0$ respectively, which is what we intended to prove.

\emph{Case} (ii): Completely similarly as in Case (i) we see that $\psi_D(L_2) \geq 0$ and $\psi_D(L_2+L_4) \geq 0$ (use the root $L_1-L_3$), and $\psi_D(L_2+L_3+L_4) \geq 0$ (use the root $L_2-L_4$). The only non-trivial case is when we want to show that $\psi_D(L_2+L_3) \geq 0$, because there is no root-hyperplane which contains both $L_1+L_3$ and $L_2+L_3$. We proceed as follows. From our assumptions we already know that $\psi_D(L_1+L_3) \geq 0$, $\psi_D(L_1+L_3)-1 <0$, $\psi_D(L_3) \geq 0$ and $\psi_D(L_1+L_2+L_3) \geq 0$. Therefore if $u(\sigma_{3,1,2})=(a,b,c,-a-b-c)$ where $a,b,c \in \Z$, then $a+c \geq 0$, $a+c-1 <0$, $c \geq 0$ and $a+b+c \geq 0$, i.e., $a=-c$, $c \geq 0$ and $b \geq 0$. But since $u(\sigma_{i,j,k})$'s form a \emph{positive} orthogonal set, we can find a non-negative integer $m$ so that $u(\sigma_{3,2,1})=m(L_2-L_1)+u(\sigma_{3,1,2})$. Thus $\left\langle u(\sigma_{3,2,1}), L_2+L_3 \right\rangle = b+c+m \geq 0$, i.e., $\psi_D(L_2+L_3) \geq 0$, and this concludes the proof.

\subsection{Proof of Theorem E for $n>4$}

We use induction to prove the theorem in the general case. Assume that the theorem is true for all $V_{SL_{n}}$ (as well as for the root-hyperplane copies of $V_{SL_{n}}$ inside $V_{SL_{n+1}}$; again see Section \ref{inductionsection}), and we want to prove it for $V_{SL_{n+1}}$.

Let $D$ be a divisor on $V_{SL_{n+1}}$ such that $\mathcal{O}(D)$ is generated by its sections. Once more we consider $D-D_{\alpha}$, where we have taken without loss of generality $\alpha = L_1-L_2$. Following the discussion in Section \ref{lemmasection}, we assume that the set (\ref{eqn:psiD}) is non-empty and that (\ref{eqn:psiDelta4}) is true. We want to prove that (\ref{eqn:psiDelta5}) holds. For a contradiction, suppose that there exists a point $v_0$ in the half-space $[\alpha \geq 0]$ such that $\psi_D(v_0)- \left\langle \alpha, v_0 \right\rangle <0 $. Since the function $\psi_D$ is continuous and piece-wise linear, we would get a contradiction if we showed that $\psi_D$ takes non-negative values on all the primitive lattice points along the rays of the fan of $V_{SL_{n+1}}$ which are contained in the half-space $[\alpha \leq 0]$. Moreover, due to the same properties of $\psi_D$, we may assume that
\begin{equation}\label{eqn:IJ}
v_0=L_1+ \sum{L_{i_j}},
\end{equation}
where $i_j$'s are distinct elements of $\{ 3,4,\ldots,n+1\}$.

We want to prove that for all the primitive lattice points of the form 
\begin{equation}\label{eqn:PQ}
L_2+\sum{L_{p_q}},
\end{equation}
where $p_q$'s are distinct elements of $\{ 3,4,\ldots,n+1\}$, we have
\begin{equation}\label{eqn:psiPQ}
\psi_D(L_2+ \sum{L_{p_q}}) \geq 0.
\end{equation}

For this, it suffices to show that there is a root $\beta$ of $\hat{G}$ such that both $v_0=L_1+\sum{L_{i_j}}$ and $L_2+\sum{L_{p_q}}$ are contained in the root hyperplane corresponding to $\beta$, because we can then apply the induction hypothesis to conclude that (\ref{eqn:psiPQ}) is true.

If there exists an $i_j$ from (\ref{eqn:IJ}) which is not equal to any of $p_q$'s appearing in (\ref{eqn:PQ}), then put $\beta = L_1-L_{i_j}$. If there exists a $p_q$ from (\ref{eqn:PQ}) which is not equal to any of $i_j$'s from (\ref{eqn:IJ}), then we put $\beta = L_2-L_{p_q}$. In both cases, (\ref{eqn:psiPQ}) follows at once.

Therefore we are only left with the possibility that, up to permutation, $i_j$'s are equal to $p_q$'s. If this is the case, then we distinguish two cases: 

Case 1). The number of $L_{i_j}$'s is greater than one, and so if, say, $L_{i_{j_1}}$ and $L_{i_{j_2}}$ appear in (\ref{eqn:IJ}) we put $\beta = L_{i_{j_1}}-L_{i_{j_2}}$ to get (\ref{eqn:psiPQ}).

Case 2). The number of $L_{i_j}$'s is less than two, and so, because $n>4$, we can find $L_{k_m}$ and $L_{k_n}$ which don't appear in (\ref{eqn:IJ}) and therefore in (\ref{eqn:PQ}) either. Put $\beta = L_{k_1}-L_{k_2}$ to see that $v_0$ and $L_2+\sum{L_{p_q}}$ belong to the root hyperplane corresponding to $\beta$. So we get (\ref{eqn:psiPQ}).

This proves Theorem E for all $n$.

\subsection{Induction process in Proof of Theorem E}\label{inductionsection}
Let us justify the induction process we used in the proof of Theorem E. Suppose $G=SL_{n+1}$ and consider $V_G$, with the initial lattice $\Z^{n+1} / \Z(1,1,\ldots,1)$ and with dual lattice $X^*(\hat{T})= \{(x_1,\ldots,x_{n+1})\in \Z^{n+1} | x_1+\ldots+x_{n+1}=0\}$. Take without loss of generality $\alpha = L_1-L_2$. Then the fan of $D_{\alpha}$ consists of the $(n-1)$-dimensional cones of the fan $\Delta$ of $V_{SL_{n+1}}$ which lie inside the hyperplane $[\alpha = 0]=\{v\in \R^{n+1}: \langle \alpha, v \rangle =0 \}.$ The dual lattice for $D_{\alpha}$ is $p_{\alpha}(X^*(\hat{T}))$. More concretely, the dual lattice for $D_{\alpha}$ is $$\{(\frac{x_2}{2},\frac{x_2}{2},x_3,\ldots,x_{n+1}):\sum_{i=2}^{n+1}x_i=0; x_i\in \Z,\forall i \},$$
because the map $p_{\alpha}$ is given by $p_{\alpha}(x_1,x_2,\ldots,x_{n+1})=(\frac{x_1+x_2}{2},\frac{x_1+x_2}{2},x_3,\ldots,x_{n+1})$. The claim is that we can treat $D_{\alpha}$ as a copy of $V_{SL_n}$ inside $V_{SL_{n+1}}$.

For ease of notation, let $L^{\vee}_{\alpha}:= p_{\alpha}(X^*(\hat{T}))$ and $L^{\vee}_n$ be the dual lattice of $V_{SL_n}$. Note that we have a bijection $\Psi: L^{\vee}_{\alpha} \rightarrow L^{\vee}_n$, given by $\Psi(x_2/2,x_2/2,x_3,\ldots,x_{n+1})=(x_2,x_3,\ldots,x_{n+1})$. Now let $\beta\neq \alpha$ be a root of $\hat{G}$. We only need to consider two cases: I) $\beta =L_2-L_3$; II) $\beta = L_3-L_4$, to prove that we can treat $D_{\alpha}$ as a copy of $V_{SL_n}$. 

I) Let $\beta =L_2-L_3$. Define $p_{\beta}$, the projection using $\beta$, in the same way as $p_{\alpha}$. We project $L^{\vee}_{\alpha}$ via $p_{\beta}$ to get the lattice $$L^{\vee}_{\alpha,\beta}:=\{ (\frac{x_3}{3},\frac{x_3}{3},\frac{x_3}{3},x_4,\ldots,x_{n+1}): \sum_{i=3}^{n+1}x_i=0; x_i\in \Z,\forall i \}.$$ If we project $L^{\vee}_n$ via the corresponding $p_{\beta}$, we get the lattice $$L^{\vee}_{n,\beta}:=\{ (\frac{x_3}{2},\frac{x_3}{2},x_4,\ldots,x_{n+1}): \sum_{i=3}^{n+1}x_i=0; x_i\in \Z,\forall i  \},$$ and we also have a natural bijection $\Psi_{\beta}:L^{\vee}_{\alpha,\beta} \rightarrow L^{\vee}_{n,\beta}$, given by $$\Psi_{\beta}(\frac{x_3}{3},\frac{x_3}{3},\frac{x_3}{3},x_4,\ldots,x_{n+1})=(\frac{x_3}{2},\frac{x_3}{2},x_4,\ldots,x_{n+1}).$$ Further, if we have an orthogonal set $(\xi_i)$ in $V_{SL_n}$, then we get an orthogonal set $(\Psi^{-1}(\xi_i))$ in $D_{\alpha}$ by applying $\Psi^{-1}$ to the former. And conversely, we apply $\Psi$ to an orthogonal set in $D_{\alpha}$ in order to get an orthogonal set in $V_{SL_n}$. Also, points in $p_{\beta}(\Conv(\xi_i)\cap L^{\vee}_{\alpha})$ are identified with points in $p_{\beta}(\Conv(\Psi^{-1}(\xi_i)\cap L^{\vee}_{\alpha})$ via the map $\Psi_{\beta}$ and, under the same map, points in $p_{\beta}(\Conv(\xi_i)) \cap p_{\beta}(L^{\vee}_{\alpha})$ are identified with those in $p_{\beta}(\Conv(\Psi^{-1}(\xi_i))) \cap p_{\beta}(L^{\vee}_n)$. 

The proof of the following Lemma then solely involves unravelling the definitions of the maps involved and the definition of a convex hull. Incidentally, it proves, as we wanted, that we can treat $D_{\alpha}$ as a copy of $V_{SL_n}$ inside $V_{SL_{n+1}}$, for which Theorem E holds.

\begin{lemma}
For the corresponding orthogonal sets $(\xi_i)$ and $(\Psi^{-1}(\xi_i))$ we have that $p_{\beta}(\Conv(\xi_i)\cap L^{\vee}_{\alpha})= p_{\beta}(\Conv(\xi_i)) \cap p_{\beta}(L^{\vee}_{\alpha})$ if and only if $p_{\beta}(\Conv(\Psi^{-1}(\xi_i)\cap L^{\vee}_{\alpha})= p_{\beta}(\Conv(\Psi^{-1}(\xi_i))) \cap p_{\beta}(L^{\vee}_n)$.
\end{lemma}

II) Let $\beta = L_3-L_4$. This is handled analogously to I) and the lemma above holds in this case, too. Let us mention that in this case we have $$L^{\vee}_{\alpha,\beta}=\{(\frac{x_3}{2},\frac{x_3}{2},\frac{x_4}{2},\frac{x_4}{2},x_5,\ldots,x_{n+1}): \sum_{i=3}^{n+1}x_i=0; x_i\in \Z,\forall i\}$$ and $$L^{\vee}_{n,\beta}=\{(x_3,\frac{x_4}{2},\frac{x_4}{2},x_5,\ldots,x_{n+1}): \sum_{i=3}^{n+1}x_i=0; x_i\in \Z,\forall i\}.$$

\subsection{Proof of Theorem D}

Let $G=GL_n$ and, as earlier, take without loss of generality $\alpha = L_1-L_2$. We have that $X^*(\hat{T})=\mathbb{Z}^n$. Suppose $D$ is a divisor on $V_G$ such that $\mathcal{O}(D)$ is globally generated and let $\{u(\sigma_i)\}_{i \in I}$ be the corresponding orthogonal set, where for each $i \in I$, $u(\sigma_i)=(a_1^{(i)},a_2^{(i)},\ldots, a_n^{(i)})\in \mathbb{Z}^{n}$ corresponds to the maximal cone $\sigma_i$ of the fan $\Delta$ of $V_G$. Then, the fact that $\{u(\sigma_i)\}_{i\in I}$ is an orthogonal set means in this case that there exists $s\in \mathbb{Z}$ such that $a_1^{(i)}+a_2^{(i)}+\ldots+ a_n^{(i)}=s,\, \forall i \in I$. 

If we denote by $P_D$ the intersection of $X^*(\hat{T})$ with the convex hull of the orthogonal set $\{u(\sigma_i)\}_{i \in I}$, then we have that (see e.g. ~\cite{fulton}, pg.66)

$$P_D=\{u\in \mathbb{Z}^n: \langle u, v\rangle \leq \langle u(\sigma_i),v\rangle , \forall v\in \sigma_i, \forall i \in I \}$$
and 
$$H^0(V_{G},\mathcal{O}(D))= \bigoplus_{u\in P_D} \C \, \chi^u.$$

Similarly, we get that

$$H^0(V_{G},i_*(\mathcal{O}(D)|_{D_{\alpha}}))= \bigoplus_{u\in P_{D_{\alpha}}} \C \, \chi^u,$$
where $P_{D_{\alpha}}= p_{\alpha}(\mathbb{Z}^n) \cap P_{\alpha}$, with $P_{\alpha}$ standing for the image of the convex hull of the orthogonal set $\{u(\sigma_i)\}_{i \in I}$ under the map $p_{\alpha}$. (Recall that the map $p_{\alpha}$ was defined in Section \ref{resultssection}).

Using the long-exact sequence (*) from Section \ref{resultssection}, Theorem D can be re-written in the following way:

$$\bigoplus_{u\in P_D} \C \, \chi^u \stackrel{\varphi}{\longrightarrow} \bigoplus_{u\in P_{D_{\alpha}}} \C \, \chi^u$$ is surjective, for all globally generated line bundles $\mathcal{O}(D)$ on $V_{G}$, where $\varphi(\chi^{u})= \chi^{p_{\alpha}(u)}$.

Let us prove Theorem D using Theorem E. Suppose $x \in P_{D_{\alpha}}$. It suffices to find a point $z$ in $P_D$ such that $x=p_{\alpha}(z)$.

First ``shift'' the orthogonal set $\{u(\sigma_i)\}_{i \in I}$ by $u(\sigma_1)$, i.e., consider the orthogonal set $\{u(\sigma_i)-u(\sigma_1)\}_{i \in I}$, which gives a new divisor $D'$ on $V_G$. Then, if we write $u(\sigma_i)-u(\sigma_1) =(b_1^{(i)},b_2^{(i)},\ldots, b_n^{(i)})$, we see that $b_1^{(i)}+b_2^{(i)}+\ldots+ b_n^{(i)}=0,\, \forall i$. So $\{u(\sigma_i)-u(\sigma_1)\}_{i \in I}$ gives an orthogonal set in $V_{SL_n}$.

Note that $p_{\alpha}(x-u(\sigma_1))$ lies in the intersection of the image of the set $$\{u\in \mathbb{R}^n: \langle u, v\rangle \leq \langle u(\sigma)-u(\sigma_1),v\rangle , \forall v\in \sigma, \forall \sigma \in \Max (\Delta) \}$$ under the map $p_{\alpha}$ and the set $p_{\alpha}(\mathbb{Z}^n)$.

Using Theorem E, we can then conclude that there exists $z'\in P_{D'}$ such that $p_{\alpha}(z')=p_{\alpha}(x-u(\sigma_1))$, where $$P_{D'}:=\{u\in \mathbb{Z}^n: \langle u, v\rangle \leq \langle u(\sigma)-u(\sigma_1),v\rangle , \forall v\in \sigma, \forall \sigma \in \Max (\Delta) \}.$$
One can easily check that $p_{\alpha}(z'+u(\sigma_1))=x$ and $z'+u(\sigma_1) \in P_D$. This means that we found $z=z'+u(\sigma_1)$ in $P_D$ such that $x=p_{\alpha}(z)$, which proves Theorem D.

\subsection{Proof of Theorem B}\label{Proof of Theorem B}

Now it is time to explain how Theorem B follows from Theorem D. We keep the same notation as in Section \ref{resultssection}, with the difference that we consider $G=GL_{n+1}$ instead of $G=GL_n$. Then $M$ must be of the form $M=GL_{n_1}\times GL_{n_2}\times \cdots \times GL_{n_r}$ where $\sum^{n+1}_{i=1}n_i=n+1$. We can also see that $\mathfrak{a}=\R^{n+1}$,
$$\mathfrak{a}_M=\{(x_1,\ldots,x_{n+1})\in \R^{n+1}:x_1=\ldots=x_{n_1}, x_{n_1+1}=\ldots = x_{n_1+n_2},\cdots,$$
$$ x_{n_1+\ldots+n_{r-1}}=\ldots=x_{n+1} \},$$
and $pr_M$ is given by averaging over each of the $r$ ``batches''. (We remind the reader that here, as before, we are considering $\mathfrak{a}_M$ as a subspace of $\mathfrak{a}$.)

The simplest case is when only one of $n_i$'s is equal to 2 and the rest are equal to 1, but we saw in Section \ref{resultssection} how Theorem B is deduced from Theorem D in that case. Now assume that only two of the $n_i$'s are equal to 2 and the rest equal 1. We can assume without loss of generality that  $n_1=n_2=2$ and $n_3=\ldots=n_r=1$, with $r=n-3$. Then $$pr_M(x_1,\ldots,x_{n+1})=(\frac{x_1+x_2}{2},\frac{x_1+x_2}{2},\frac{x_3+x_4}{2},\frac{x_3+x_4}{2},x_5,\ldots,x_{n+1})$$
and so $pr_M$ is just the composition $p_{L_3-L_4}\circ p_{L_1-L_2}$. Therefore when we apply Theorem D to $D_{L_1-L_2}$, where we consider $D_{L_1-L_2}$ as a copy of $V_{GL_n}$ and $\alpha=L_3-L_4$, we find that Theorem B follows in this case, too.

The next case would be to consider, without loss of generality, $n_1=3$ and $n_2=\ldots=n_r=1$, where $r=n-2$. Now we see that $pr_M$ is the composition $p_{L_2-L_3}\circ p_{L_1-L_2}$ and we get our desired conclusion again using Theorem D.

In general we write $pr_M$ as a composition of a finite number of $p_{\alpha}$'s for distinct roots $\alpha$ and then apply Theorem D as many times as there are roots $\alpha$ appearing in that composition. If we want to be very specific, then, for $M=GL_{n_1}\times GL_{n_2}\times \cdots \times GL_{n_r}$ we can see that $pr_M$ agrees with the composition of maps $q_1,q_2,\ldots,q_r$ where $$q_k=p_{L_{(n_k-1)}-L_{n_k}}\circ p_{L_{(n_k-2)}-L_{(n_k-1)}} \circ \ldots \circ p_{L_{(n_{k-1}+1)}-L_{(n_{k-1}+2)}},$$
with the convention $n_0=0$. This concludes the proof of Theorem B.

\section{Proof of Theorem C}

Consider $V_G$ for $G={G_2}$. Just as in the case of $SL_n$, we would like to describe the fan $\Delta$ of $V$. But first, the initial lattice is $L=\Z^3 / \Z (1,1,1)$ and so the dual one is $L^{\vee} = \{(x_1,x_2,x_3) \in \Z^3:x_1+x_2+x_3=0 \}$. Identify the rays in $\Delta$ with their corresponding minimal lattice points. Then $\Delta (1)$ consists of twelve points: $v_1=(1,0,0), v_2=(1,0,-1), v_3=(0,0,-1), v_4=(0,1,-1), v_5=(0,1,0), v_6=(-1,1,0)$ and $v_{i+6}=-v_i$, for $i=1,\ldots,6$. There are twelve maximal cones: $\sigma_i$, $i=1,2,\ldots,12$, where the rays of $\sigma_i$ are $v_i$ and $v_{i+1}$ (with $v_{13}:=v_1$). Let $D$ be a $\hat{T}$-Cartier divisor on $V$ corresponding to a Weyl orbit (see Section \ref{resultssection} for the precise definition) and $\alpha$ a root of $\hat{G}$. For every maximal cone $\sigma_i$ denote by $u(\sigma_i)$ the character corresponding to the divisor $D$. Obviously, $u(\sigma_i)$'s form a positive orthogonal set. We want to prove that $$H^1(V_G,\mathcal{J}_{D_{\alpha}} \otimes \mathcal{O}(D))=0.$$
This is true if the orthogonal set corresponding to $D$ is \emph{strictly positive}, i.e., if $u(\sigma_i)$'s are distinct and they form a positive orthogonal set. Indeed, use Lemma \ref{Dalpha2} and the definition of orthogonal sets to conclude that the orthogonal set corresponding to $D-D_{\alpha}$ is still positive (but not necessarily strictly positive), i.e., $\mathcal{J}_{D_{\alpha}} \otimes \mathcal{O}(D)$ is globally generated. It is also obvious that if all $u(\sigma_i)$'s are equal to each other, then we get a degenerate case for which the result is true as well.

Therefore the only problem might arise if only some but not all of the $u(\sigma_i)$'s are the same. Using symmetries of the current root system, there are only four cases we need to consider, two for a short root and two for a long root:
\begin{itemize}
\item [(a)] $\alpha = (1,-1,0)$ and there exists a natural number $n$ so that $$u(\sigma_2)=u(\sigma_3)=-u(\sigma_8)=-u(\sigma_9) =n(1,1,-2),$$ $$u(\sigma_4)=u(\sigma_5)=-u(\sigma_{10})=-u(\sigma_{11}) =n(-1,2,-1),$$
$$u(\sigma_6)=u(\sigma_7)=-u(\sigma_{12})=-u(\sigma_1) =n(-2,1,1);$$

\item [(b)] $\alpha = (2,-1,-1)$ and there exists a natural number $n$ so that $$u(\sigma_2)=u(\sigma_3)=-u(\sigma_8)=-u(\sigma_9) =n(1,1,-2),$$ $$u(\sigma_4)=u(\sigma_5)=-u(\sigma_{10})=-u(\sigma_{11}) =n(-1,2,-1),$$
$$u(\sigma_6)=u(\sigma_7)=-u(\sigma_{12})=-u(\sigma_1) =n(-2,1,1);$$

\item [(c)] $\alpha = (1,-1,0)$ and there exists a natural number $n$ so that $$u(\sigma_1)=u(\sigma_2)=-u(\sigma_7)=-u(\sigma_8) =n(1,0,-1),$$ $$u(\sigma_3)=u(\sigma_4)=-u(\sigma_9)=-u(\sigma_{10}) =n(0,1,-1),$$
$$u(\sigma_5)=u(\sigma_6)=-u(\sigma_{11})=-u(\sigma_{12}) =n(-1,1,0);$$

\item [(d)] $\alpha = (2,-1,-1)$ and there exists a natural number $n$ so that $$u(\sigma_1)=u(\sigma_2)=-u(\sigma_7)=-u(\sigma_8) =n(1,0,-1),$$ $$u(\sigma_3)=u(\sigma_4)=-u(\sigma_9)=-u(\sigma_{10}) =n(0,1,-1),$$
$$u(\sigma_5)=u(\sigma_6)=-u(\sigma_{11})=-u(\sigma_{12}) =n(-1,1,0).$$
\end{itemize}

Let us mention again that we want to prove that $H^1(V_G,\mathcal{J}_{D_{\alpha}} \otimes \mathcal{O}(D))=0$. As we explained in Section \ref{resultssection}, we have a long-exact sequence
$$... \longrightarrow H^0(V_G, \mathcal{O(D)}) \stackrel{\varphi}{\longrightarrow} H^0(V_G,i_*(\mathcal{O(D)}|_{D_{\alpha}})) \longrightarrow H^1(V_G,\mathcal{J}_{D_{\alpha}} \otimes \mathcal{O(D)}) \longrightarrow 0,$$
where terms on the right vanish because $\mathcal{O(D)}$ is generated by its sections. Therefore it suffices to prove that $\varphi$ is surjective in each of the cases above.

Case (a): We have that (cf. ~\cite{fulton}, pg.66) $$H^0(V_G, \mathcal{O(D)})= \bigoplus_{u\in P_D}\C \chi^u,$$ where $P_D=\{u\in L^{\vee}: \langle u, v\rangle \leq \langle u(\sigma_i),v\rangle , \forall v\in \sigma_i, \forall i=1,\ldots,12 \}$.
Also, $$H^0(V_G,i_*(\mathcal{O(D)}|_{D_{\alpha}})= \bigoplus_{u\in P_{D_{\alpha}}}\C \chi^u,$$ where $P_{D_{\alpha}}=\{(\frac{k}{2},\frac{k}{2},-k): k\in [-2n,2n]\cap \Z \}$. The map $\varphi$ is given by $\chi^{(a,b,-a-b)} \mapsto \chi^{(\frac{a+b}{2},\frac{a+b}{2},-a-b)}$, where $u=(a,b,-a-b) \in P_D$. It is clear that if $k$ is even, we have $\varphi(\chi^{(\frac{k}{2},\frac{k}{2},-k)})=\chi^{(\frac{k}{2},\frac{k}{2},-k)}$ and $(\frac{k}{2},\frac{k}{2},-k) \in P_D$. If $k$ is odd, then we have  $\varphi(\chi^{(\frac{k+1}{2},\frac{k-1}{2},-k)})=\chi^{(\frac{k}{2},\frac{k}{2},-k)}$. Therefore to prove that $\varphi$ is surjective, it is enough to show that $(\frac{k+1}{2},\frac{k-1}{2},-k) \in P_D$ for all odd integers $k$ between $-2n$ and $2n$. But this is easily seen to be true by checking that $\langle (\frac{k+1}{2},\frac{k-1}{2},-k), v_j\rangle \leq \langle u(\sigma_i),v_j\rangle$, where $v_j$'s are the respective primitive lattice points along the rays of our fan $\Delta$, written explicitly at the beginning of this section.

Case (b): Keeping the same notation as in Case (a) we now have that $P_{D_{\alpha}}= \{ (0, \frac{k}{2},-\frac{k}{2}): k\in [-3n,3n] \cap \Z \}$. $\varphi$ is given by $\chi^{(a,b,-a-b)} \mapsto \chi^{(0, \frac{a}{2}+b,-\frac{a}{2}-b)}$, where $(a,b,-a-b)\in P_D$. If $k$ is even then $(0,\frac{k}{2},-\frac{k}{2}) \in P_D$ and $\varphi$ fixes $\chi^{(0,\frac{k}{2},-\frac{k}{2})}$. If $k$ is odd then we have that $\varphi (\chi^{(1,\frac{k-1}{2},-\frac{k+1}{2})})=\chi^{(0, \frac{k}{2}, -\frac{k}{2})}.$ It is an easy calculation to see that for all odd $k$ between $-3n$ and $3n$ we get $(1,\frac{k-1}{2},-\frac{k+1}{2}) \in P_D$ and therefore $\varphi$ is surjective, which we wanted to prove. 

Case (c): This time we have $P_{D_{\alpha}}= \{ (\frac{k}{2},\frac{k}{2},-k): k\in [-n,n] \cap \Z \}$ and $\varphi$ is the same as in case (a). Clearly $(\frac{k}{2},\frac{k}{2},-k)$ is in $P_D$ and $\chi^{(\frac{k}{2},\frac{k}{2},-k)}$ is fixed by $\varphi$, whenever $k$ is even. If $k$ is odd then consider $(\frac{k+1}{2},\frac{k-1}{2},-k)$, which is easily seen to lie in $P_D$. $\varphi$ is surjective because it maps $\chi^{(\frac{k+1}{2},\frac{k-1}{2},-k)}$ to $\chi^{(\frac{k}{2},\frac{k}{2},-k)}$.

Case (d): In our final case, we have $P_{D_{\alpha}}= \{ (0,\frac{k}{2},-\frac{k}{2}): k\in [-2n,2n] \cap \Z \}$. $\varphi$ is the same as in Case (b) and it fixes $\chi^{(0,\frac{k}{2},-\frac{k}{2})}$ for $k$ even. Clearly $(0,\frac{k}{2},-\frac{k}{2}) \in P_D$ and it is not difficult to see that if $k$ is odd then $(1, \frac{k-1}{2},-\frac{k+1}{2}) \in P_D$ as well. Moreover, still for $k$ odd, $\varphi (\chi^{(1, \frac{k-1}{2},-\frac{k+1}{2})}) =  \chi^{(0,\frac{k}{2},-\frac{k}{2})}.$ Hence $\varphi$ is surjective.

This concludes the proof of Theorem C and therefore, as explained in Section \ref{resultssection}, of Theorem A, too.

\subsection{}
We would like to note that the proof of Theorem D was different from that of Theorem C because in the former we used \emph{all} positive orthogonal sets whereas in the latter we only used some and not all positive orthogonal sets, namely only the Weyl-orbit ones. More concretely, in Section 3,  after Lemma \ref{importantlemma}, we saw that the vanishing of the first cohomology group there was implied by the vanishing of its corresponding 0-eigenspace (under the action of our torus). 
This is not true for $G=G_2$. The crucial step was to note that by ``shifting" a positive orthogonal set, say $\{u(\sigma)\}$, by a character, say, $u$, we get a positive orthogonal set $\{u(\sigma)-u\}$; however, it is clear that a shift of a Weyl-orbit is in general not a Weyl-orbit.

\subsection{}
At the very end, we give an example to show that Theorem B fails in the case of $G_2$. Consider $V_{G_2}$ and let $\alpha=(1,-1,0)$. Let $D$ be the positive orthogonal set (i.e., the globally generated torus-equivariant divisor on $V_{G_2}$) given by $u(\sigma_i)=(1,0,-1)$ for $i \in \{1,2,3,4,5,12 \}$ and $u(\sigma_i)=(0,-1,1)$ for $i \in \{6,7,8,9,10,11 \}$. Following the same notation as before, we see that in our case $$ P_{D} = \{ (1,0,-1), (0,-1,1) \} ,$$ and $$ P_{D_{\alpha}} = \{ (\frac{1}{2}, \frac{1}{2}, -1), (0,0,0), (-\frac{1}{2}, -\frac{1}{2}, 1) \} .$$ Therefore $\Dim H^0(V_{G_2}, \mathcal{O(D)}) =2$, but $\Dim H^0(V_{G_2},i_*(\mathcal{O(D)}|_{D_{\alpha}})) =3$. Hence the map $\varphi : H^0(V_{G_2}, \mathcal{O(D)}) \rightarrow H^0(V_{G_2},i_*(\mathcal{O(D)}|_{D_{\alpha}}))$, induced by the projection $p_{\alpha}$ along $\alpha$, is not surjective, and Theorem B is no longer true if we put $G_2$ instead of $GL_n$.

\end{document}